\documentclass[10pt]{article}
\usepackage[english]{babel}
\usepackage{latexsym}
\usepackage[latin1]{inputenc}
\usepackage{amsfonts}
\usepackage{mathrsfs}
\usepackage{eucal}
\usepackage{exscale}
\usepackage{makeidx}

\usepackage{graphicx}
\usepackage[T1]{fontenc}
\usepackage{amsmath,amssymb}

\usepackage[all]{xy}
\usepackage{pgf,tikz}
\usepackage{mathrsfs}
\usetikzlibrary{arrows}

\newtheorem{thm}{Theorem}
\newtheorem{rem}{Remark}
\newtheorem{rems}{Remarks}
\newtheorem{df}{Definition}[section]

\newtheorem{cor}{Corollary}
\newtheorem{lem}{Lemma}
\newtheorem{fait}{Fact}[section]

\newtheorem{exemples}{Examples}

\newcommand{\preuve}{\indent {\it Proof:}\hspace{4mm}}

\title{Positive amalgamation}
\author{Mohammed Belkasmi}

\begin{document}
\maketitle
\abstract 
 
We  study the amalgamation property in positive logic. We give  some connections
 between the amalgamation 
property and Robinson theories, model-complete 
theories and the Hausdorff property.
\vspace{0.5cm}

\section*{Introduction}
The  positive model theory is the 
study of the h-inductive theories through special classes of 
 of their models, particularly the class of existentially closed models,
 the class of h-maximal models, and the class of amalgamation basis.\\ 
The notion of positive model theory was first introduced by 
Ben Yaacov following the line of research
on inductive classes of structures accomplished by Pillay.

In this paper we continue the study of the amalgamation property
 begun in \cite{ana, begnacpoizat}. 
In the first section we recall the basic concepts of positive logic and we  prove some properties of the classes of positive 
closed and h-maximal models. In the second section we give
some results
about the amalgamation  bases and amalgamation class and their connections with Robinson theories and Hausdorff property.
The third section is devoted to the axiomatization of the amalgamation 
basis and the amalgamation property.

\section{Positive logic}
In the first part of this section we will recall various definitions and notions of positive logic. For further detail, \cite{begnacpoizat} is sufficiently complete reference.

Consider $L$ a first-order language. A  formula is said to be positive 
 when it is obtained from the atomic formulas by the use of 
$\exists, \vee, \wedge$. It can be written in prenex form as 
$\exists\bar x\ \phi(\bar x, \bar y)$, where $\phi$ is
positive quantifier-free. In this case the variables $\bar y$ are said 
to be free. 

A sentence is a first order formula without free variables.
A sentence is said to be $h$-universal
if it is the negation of a positive sentence. The set of h-universal 
sentences is closed under conjunction and disjunction.\\
A sentence is said to be $h$-inductive if it is a finite conjunction
of sentences of the form
\[
\forall\bar x(\exists\bar y \varphi(\bar x,\bar y)
\rightarrow\exists z \psi(\bar x,\bar z)) \ ,
\]
where $\varphi$ and $\psi$ are 
quantifier-free and positive formulas.

Let $A$ and $B$ be  $L$-structures. A mapping $f$
from $A$ into $B$ is a homomorphism if and only if for evry tuple
$\bar a$ from $A$,
and for every atomic formula $\phi$; $A\models\phi(\bar a)$
implies $B\models\phi(f(\bar a))$.
In such a case $B$ is said to be a continuation of $A$.\\
A homomorphism $f$ is an embedding  if and only if for every tuple
 $\bar a$ from $A$,
$\bar a$ and $f(\bar a)$ satisfy the same atomic formulas.\\
A homomorphism is said to be an immersion whenever $\bar a$
and $f(\bar a)$ satisfy the same positive formulas.

A class of L-structures is said to be $h$-inductive if it is closed
with respect to inductive limits of homomorphisms. 
A  theory $T$ is said to be $h$-inductive  
 if  the class of models of $T$ is $h$-inductive.

\begin{fait}[{\cite[Th\'eor\`eme 23]{begnacpoizat}}]\label{ind}
A theory is $h$-inductive
if and only if it is axiomatized by $h$-inductive sentences.
\end{fait}

\subsection{Positively existentially closed models}
In this  subsection, we will give some
properties and characterizations of special
classes of structures of an h-inductive theory. Namely
the classes of positively  existentially closed
 models  (in short, pc), and h-maximal models. 
\begin{df}
A member $M$ of a class $\Gamma$ of $L$-structures
is said to be pc (respectively h-maximal) in $\Gamma$,
 if every homomorphism from $M$ into 
a member
of $\Gamma$ is an immersion (respectively an embedding).
\end{df}
In \cite{begnacpoizat} it is shown that the  pc models exist for 
any consistent h-inductive theory. We have the 
useful following fact.
\begin{fait}[{\cite[Th\'eor\`eme 1, lemme 12] {begnacpoizat}}]\label{pecconti}:
\begin{itemize}
\item Every member of an $h$-inductive class is continued 
in a pc member of the class.
\item The class of pc models of an $h$-inductive theory $T$ is
 $h$-inductive.
\end{itemize}
\end{fait}

\begin{df}[{\cite{begnacpoizat}}]
Two $h$-inductive theories are said to be companions
if they have the same $h$-universal consequences.
\end{df}
\begin{fait}[{\cite[lemme 7]{begnacpoizat}}]\label{comp}
Two $h$-inductive theories are companions if and only if
they have the same pc models.
\end{fait}
\begin{rem} 
 Two companion theories do not necessarily share the same
h-maximal models.
\end{rem}

The  $h$-inductive
theory $T$ has a maximal companion denoted $T_k(T)$ called
the Kaiser's hull of $T$; it
 is the set of all $h$-inductive sentences
true in every pc model  of $T$. Similarly,
$T$ has a minimal companion denoted $T_u(T)$, formed by its
$h$-universal consequences.\\
Let $A$ be a L-structure. We denote by $T_i(A)$ the set of h-inductive
sentences which are true in $A$. \\
Let $L'$ be the language  obtained from $L$ by  adding the elements 
of $A$ as constants to $L$.
We denote by 
$T_k(A)$ the set of h-inductive sentences in the language $L'$ which 
are true in $A$.
\begin{df}
An $h$-inductive theory $T$ is said to be model-complete
if all models of $T_k$ are pc models of $T$. In other words, the class
of pc models is axiomatised by the Kaiser's hull $T_k$.\\
An h-inductive class $\Gamma$ is said to be complete if it has the 
jhp property (ie for every $A, B$ elements of $\Gamma$ there is $C$
an element of  $\Gamma$ in which $A$ and $B$ are continued)
\end{df}
\begin{lem}
An h-inductive  theory $T$ is complete if and only if 
its pc models  have the same h-inductive theory.

\end{lem}
\preuve 
Let $A$ and  $B$  be two  pc models of $T$, since $T$ is a complete theory. $A$ and $B$ are immersed in a pc model $C$ of $T$. Thereby 
$A$, $B$ and $C$ satisfy the same h-inductive theory.

For the other direction, since any two pc models of $T$ share the same 
h-inductive theory, by compactness they have a common continuation in the class of models of $T$.
  

\begin{lem}\label{finitecomplete}
Let $T$ be a complete h-inductive theory which has
a finite pc model $A_e$. Then $T$ is  model-complete 
and the class of its pc models is reduced to $A_e$ up to isomorphism.
\end{lem}
\preuve Let $B$ be 
a pc model of $T$.  We claim that  $\mid A_e\mid=\mid B\mid$.
Indeed, since $T$ is complete,  then its pc models  
 have the same h-inductive theory $T_k(T)$.
 On the other hand,
since $A_e$  satisfies $\chi$  the following h-inductive sentence 
$$\forall x_1,\cdots,x_{\mid A\mid+1}\, 
(\bigvee_{1\leq i\neq j\leq \mid A\mid+1}x_i=x_j)$$
Then $B\vdash\chi$, thereby $|B|\leq |A_e|$. 
Similarly, we have  $|A_e|\leq |B|$.\\
Moreover, since  $A_e$ and $B$ are immersed in a
pc model $C$ which has the same cardinality of $A_e$. 
Then $A_e$ and 
$B$ are isomorph.

\exemples\label{exemple1} :
\begin{itemize}
\item Consider the language $L=\{R_i : \, i<\omega\}$, 
where $R_i$ are 1-ary relations symbol. 
Let $T_\omega$
be the h-universal theory  
$$\{\neg\exists x\ (R_i(x)\wedge R_j(x))\, \mid \, i,j<\omega, i\neq j\}$$
The L-structure $A_e=\{x_i\,\mid\, i<\omega\}$ such that  for all
$i<\omega$ we have $A_e\models R_i(x_i)$ is the  unique pc model of
 $T_\omega$. 
Then $T_\omega$ is not model-complete theory.\\
The h-maximal models of $T_\omega$ 
are the  substructures of  $A_e$.

\item  Let $L$ be the relational language  $\{S(x, y)\}$. Let $T_3$ 
be the h-inductive theory
$$\{\forall\  x_1,\cdots,x_5\, (\bigwedge_{i=2}^5 S(x_1, x_i)\longrightarrow
 \bigvee_{1\leq i\neq j\leq 5} x_i=x_j),\ \  \neg\exists x\, 
 S(x,x)\}$$
$T_3$ is a complete theory. The $L$-structure
 $A_e=\{a, b, c, d\}$ such that 
  $$A_e\vdash \forall xy\ ((x=y)\vee S(x, y))$$
 Is a pc model of $T_3$.  By the lemma 
 \ref{finitecomplete}, $A_e$  is the unique pc model of $T_3$. Thereby
 $T_3$ is model-complete.\\
The h-maximal models of $T_3$ are the substructures of $A_e$.
\item This example is given in \cite{almaz}.
 Let $L$ be the functional language $\{f\}$. 
 Let $T_s$ be the h-inductive theory
$$T_s=\{\neg\exists x\ (f(x)=x) \}$$
 The unique pc model of 
$T_s$ is the $L$-structure $A_e$ formed by all p-cycles, where 
$p$ is a prime number. Thereby $T_s$ is not model-complete theory. 
The h-maximal models of $T$ are the substructures of $A_e$.

\end{itemize}
Generally the class of h-maximal models of $T$  is not 
always the class of substructures of the pc models. A simple instance of this situation is given  by the theory
of rings. The domain  $\mathbb{Z}$ is a substructure of the pc
model $\mathbb{C}$, but it is not a h-maximal model of rings theory. 
\subsection{Types space}
\begin{df}
Let $T$ be an $h$-inductive theory in a language $L$.
An $n$-type is a maximal set of positive formulas
in $n$ variables that is consistent with $T$.
\end{df}

A model $A$ of an h-inductive theory $T$ is pc model, if and only if for every $\bar a\in A$,
the set of positive formulas satisfied by $\bar a$
is a type.\\
Therefore, if   a pc model $A$ of  a theory $T$ does not satisfy
 $\phi(\bar a)$,
where $\phi$ is a positive formula and $\bar a\in A$. Then there is
$\psi$ a positive formula such that $A\models\psi(\bar a)$
and $$T\vdash\neg\exists x(\phi(\bar x)\wedge\psi(\bar x)).$$
Let $\phi$ be a positive formula.
We denote by $Res_T(\phi)$, the resultant of $\phi$ modulo $T$, is the
set of positive formulas $\psi$ such that 
$$T\vdash\neg\exists x(\phi(\bar x)\wedge\psi(\bar x)).$$

Let $A$ be a model of an $h$-inductive theory $T$
and $\bar a\in A$. We denote by $F_A(\bar a)$ (resp $D_A(\bar a)$) the
set of  positive formulas (resp quantifier-free positive formulas)
 satisfied by $\bar a$ in $A$.
We denote by $S_n(T)$ the space
of $n$-types of a theory $T$. 

One defines a topology on the space $S_n(T)$  for which
the basis of closed sets are  the sets
$F_\varphi$, where $\varphi$ ranges over the entire
set of positive formulas, and
\[
F_\varphi\ =\ \{\ p\in S_n(T)\ |\ p\vdash\varphi\ \}\ 
\]

\begin{df}
An h-inductive theory $T$ is said to be 
Hausdorff if and only if for every natural number n,
 the topological space $S_n(T)$ is Hausdorff.
 A L-structure $M$ is said to be Hausdorff if and only if the theory 
 $T_k(M)$ is Hausdorff  in the language obtained from $L$
 by adding the elements of $M$ as constants.
 \end{df}
 
Note that,  The topological spaces  $S_n(T)$ are quasi-compact,
but generally  are not  Hausdorff.
\begin{df}
Let $(\varphi, \psi)$ be a pair of positive formulas such that 
$T\nvdash \forall\bar x \varphi(\bar{x})
\vee\forall\bar{y}\psi(\bar y)$.\\
 $(\varphi, \psi)$ is said to be $T$-complementary if and only if
$T\vdash\neg\exists\bar x\ (\varphi(\bar x)\wedge\psi(\bar x))$ and 
$T\vdash\forall\bar{x}\ (\varphi(\bar{x})\vee\psi(\bar{x}))$.
In this case each of the both 
formulas is said to be $T$-complemented.\\
 We say that $(\varphi, \psi)$ is $T$-weakly complementary if and only if
$T\vdash\forall\bar{x}\ (\varphi(\bar{x})\vee\psi(\bar{x}))$. In this case each of the both 
formulas is said to be $T$-weakly complemented.
\end{df}

Consider an h-inductive theory $T$. For every positive formula 
$\varphi$. We denote by $O_\varphi$ the set  of types $p\in S(T)$ such that
$p\nvdash\varphi$. We have the following properties:\\
For every positive formula $\varphi$;
 $$O_\varphi=(F_\varphi)^c=\bigcup_{\phi\in Res_T(\varphi)}F_\phi.$$
A pair of positive formulas $(\varphi, \psi)$ is $T_k(T)$-complementary
if and only if $O_\varphi= F_\psi$.\\
A theory $T$ is model-complete if and only if every positive 
formula is $T_k(T)$-complemented.

\begin{lem}
Let $T$ be an h-inductive theory. $T_k(T)$ is model-complete if and only if
for every positive formula $\phi$, there is a
 $T_k(T)$-complementary pair of 
positive formulas $(\varphi, \varphi^c)$  such that:
$$T_k(T)\vdash\forall\bar{x}\ (\phi(\bar{x})
\rightarrow \varphi(\bar{x}))$$ 
and for every $\phi'\in Res_T(\phi)$ we have 
$$T_k(T)\vdash\forall\bar{x}\ (\phi'(\bar{x})
\rightarrow \varphi^c(\bar{x})).$$
\end{lem}
\preuve 
Suppose that $T_k(T)$ is model-complete. Then each positive 
formula  $\phi$ has a $T_k(T)$-complement $\phi^c$.
 We take $\varphi=\phi$ and $\psi=\phi^c$. 
 The considered formulas  satisfy the desired properties.
 
 For the reverse direction. Let $\phi$ be a  positive 
 formula. By hypothesis  there is a $T_k$-complementary
  pair $(\varphi, \psi)$ of 
  positive formulas 
  such that; for every $\phi'\in Res_T(\phi)$, the theory $T_k(T)$
 satisfies the following h-inductive sentences;
 $$\forall\bar{x}\ (\phi(\bar{x})
\rightarrow \varphi(\bar{x})); \ \ \  
\forall\bar{x}\ (\phi'(\bar{x})
\rightarrow \varphi^c(\bar{x}))$$
This means that in the topological space $S(T)$ we have 
$F_{\varphi^c}=O_{\phi}$. Then $\phi$ is $T_k(T)$-complemented. 
Thereby  $T$ is model-complete. 
\subsection*{ Positive Robinson theories}
 
 \begin{df}\label{posrobdefn} Let $T$ be an $h$-inductive
theory.  $T$ is said to be positive Robinson if
it satisfies the following condition:\\
For  any pc models $A$ and $B$ of $T$;  

 $$\forall (\bar a, \bar b)\in A\times B;\ \ \ 
 D_A(\bar a)\subseteq D_B(\bar b)
 \Rightarrow tp(\bar a)=tp(\bar b).$$
 \end{df}
\textbf{Remark}:\\
If $T$ is positive Robinson theory, then all
its companion theories are of Robinson.\\ 
\begin{exemples} Let 
$T_\omega, T_3$ and $T_s$ be  the theories given in the example 
\ref{exemple1}; 
\begin{itemize}
\item $T_\omega$, $T_3$ and $T_s$ are positive Robinson theories.
\item The theory of commutative rings is a positive Robinson theory.
\item The theory of groups is a positive Robinson theory. Indeed
in the framework of  positive logic,
 the unique pc model of this theory is the trivial group $\{e\}$.
\end{itemize}
\end{exemples}

\begin{lem}\label{robsubpec}
Let $T$ be a positive Robinson theory. The h-maximal 
models of $T_u$ are the substructures of the pc models 
of $T$.
\end{lem} 
\preuve
Consider the following  the diagram 
\[
\xymatrix{
    A \ar[r]^{e} \ar[d]_{f} & A_e  \\
    B \ar[r]_{g} & B_e
  }
\]

Where $A_e$ is a pc  model of $T$,   $A$  a substructure of 
$A_e$, $B$ a model of $T_u$,  
$B_e$  a pc model of $T$ in which $B$ is 
continued by a 
homomorphism $g$, and 
$f$ a  homomorphism from $A$ into $B$.

  Let $\bar a$ be a tuple of $A$. For every  quantifier-free
positive formula $\varphi$ such that $B\models\varphi(f(\bar a))$ we have 
$B_e\models \varphi(g\circ f(\bar a))$. Thereby; 
$$D_{A_e}(\bar a)=D_A(\bar a)\subset D_{B_e}(g\circ f(\bar a))$$ 
Moreover, since $T$ is a positive Robinson theory we obtain 
$D_A(\bar a)=D_{B_e}(g\circ f(\bar a))$ for each $\bar a\in A$, 
this means that  $g\circ f$
is an embedding, then  $f$ is an embedding. Consequently $A$ is 
a h-maximal model of $T_u$.

Reciprocally,  since each model $A$ of $T_u$ is continued into a pc model
of $T_u$, Then every h-maximal model of $T_u$  is embedded 
in some pc model of $T$.

\section{Amalgamation bases}
\begin{df}\label{defab}
Let $\Gamma$ be a class of L-structures. An element $A$ of $\Gamma$ 
is said to be  an amalgamation basis of $\Gamma$
 if and only if for all $B, C$ in $\Gamma$, and $f, g$ 
homomorphisms 
from $A$ respectively into $B$ and $C$, there exist 
$D\in\Gamma$ and $f', g'$ homomorphisms 
 such that
the following diagram commutes;
\[
\xymatrix{
    A \ar[r]^{f} \ar[d]_{g} & {B} \ar[d]^{g'} \\
    C \ar[r]_{f'} & {D}
  }
\]
 We say that $\Gamma$ has the amalgamation 
property if every element of $\Gamma$ is an amalgamation basis. 
\end{df}
\begin{rems}
 Let $T, T'$ be two h-inductive companion theories.
\begin{itemize}
\item. If $T$ has  the 
amalgamation property we can not say the same for $T'$, except 
if $T\subset T'$.  
\item 
The classes of amalgamation bases of $T$ and $T'$
are not necessarily equal.
\end{itemize} 
\end{rems}
\begin{exemples}
\begin{itemize}
\item Let $A$ be a L-structure. $A$ is an amalgamation 
basis of $T_u(A)$, where $T_u(A)$ is the h-universal theory of
$A$ in the language $L'$ obtained from $L$ by adding the elements of 
$A$  as constants. Indeed, since $A$ is immersed in each model 
of $T_u(A)$, the fact that $A$ is an amalgamation basis of 
$T_u(A)$ follows from the fact \ref{amalgamationassy}.  
\item Every pc model of an h-inductive theory $T$ 
is an amalgamation basis 
of $T$.
\item Let $T$ be a model-complete theory. Since each model of 
$T_k(T)$ is a pc model of $T$, 
 then $T_k(T)$ has the amalgamation property.
 \item \label{2relations} 
Let $P$ and $Q$  two   1-ary relation symbols. Let $L$ be 
the relational language $ \{P, Q\}$. Let $T$ be 
 the h-universal theory  consisting  in the axiom 
  $\neg\exists xy (P(x)\wedge Q(y))$.\\
$T$  has exactly two pc models
$A_e=\{a\}$ and  $B_e=\{b\}$ such that, 
$A_e\models  P(a)$ and $B_e\models Q(b)$. 
Thereby $T$ is model-complete. 
 
$T$ has the amalgamation property. Indeed,  
let $A, B, C$ be  models of $T$. Suppose that 
$A$ is  continued  in $B$ and  $C$, and $A\vdash\exists xP(x)$.
Then $B$ and $C$  satisfy $\exists xP(x)$, 
thereby the following diagram
commutes:
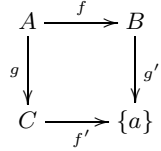
\[
\xymatrix{
    A \ar[r]^{f} \ar[d]_{g} & {B} \ar[d]^{g'} \\
    C \ar[r]_{f'} & {\{a\}}
  }
\]
Where $\{a\}$ is the pc model $A_e$. 
\item $T_\omega$ and $T_k(T_\omega)$ the Kaiser's hull of $T_\omega$
do not  have  the amalgamation property. 
\item Let $T_n$ be a theory consisting in a finite  subset of 
the set of h-universal sentences of $T_\omega$. Since $T_n$  
is complete and has a finite pc models, by the lemma \ref{finitecomplete},
$T_n$ is model- complete. Thereby $T_k(T_n)$  has the amalgamation property.
\item $T_3$ does not have the amalgamation property. Indeed, let $A, B, C$ be  the models of  $T_3$ given in the figure below:\\
$
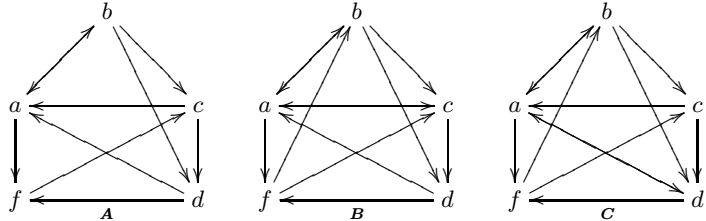
\xymatrix{
    &b \ar[rd] \ar[ld]\ar[ddr] &   \\
    a\ar[ru]\ar[d]& & c\ar[d]\ar[ll] \\
    f\ar[rru]&&d\ar[llu]\ar[ll]^{\textbf{A}}
  } \hspace{0.5cm}
  \xymatrix{
    &b \ar[rd] \ar[ld]\ar[ddr] &   \\
    a\ar[ru]\ar[d]\ar[rr]& & c\ar[d]\ar[ll] \\
    f\ar[rru]\ar[ruu]&&d\ar[llu]\ar[ll]^{\textbf{B}}
  }
  \hspace{0.5cm}
 \xymatrix{
    &b \ar[rd] \ar[ld]\ar[ddr] &   \\
    a\ar[ru]\ar[d]\ar[rrd]& & c\ar[d]\ar[ll] \\
    f\ar[rru]\ar[ruu]&&d\ar[llu]\ar[ll]^{\textbf{C}}
  } 
$

$A$ is continued into  $B$ and $C$, and 
the diagram resulting can  not be amalgamate. Note that 
$T_k(T_3)$ has the amalgamation property.
\item $T_s$ does not have the amalgamation property. Indeed,   
 $\mathbb{Z}$ is continued  onto two distinct p-cycles, but the diagram resulting can not be amalgamate   by a model of $T_s$.
\item The fields constitutes an h-inductive  subclass of amalgamation bases in the class of commutative  rings.
\end{itemize}
\end{exemples}
\begin{lem}\label{lem amal basis rings}
If $A$ is an amalgamation basis of the theory of 
rings then $A$ is simple.
\end{lem}
\preuve Suppose that $A$ is a not-simple  
 amalgamation basis of rings. Let 
$I$ be a maximal proper ideal of $A$, $x\in A-I$ and 
let $(x)$ be  the ideal generated by $x$.\\
Since $A$ is an amalgamation basis and continued into $A/I$ and 
$A/(x)$;
there is $B$ a ring such that the following diagram commutes 
 \[
\xymatrix{
    A \ar[r]^{f} \ar[d]_{g} & {A/(x)} \ar[d]^{f'} \\
    B/I\ar[r]^{e}& B
  }
\]
Where $f,f'$ and $g$ are homomorphisms,  $e$ an embedding.\\ 
Consequently, we have $e(g(x))=f'(f(x))=0$, this   implies that 
$x\in I$. Contradiction, thereby $A$ is a simple ring.

\begin{rem}
 $B$ is an amalgamation basis of the theory of commutative 
 rings if and only if it is a field.
\end{rem}

The following facts give a useful characterisation  of positive
Robinson theories and amalgamation bases. The fact \ref{amalgamationassy} is a modified version of the so-called 
"asymetric amalgamation" given in \cite{begnacpoizat}. 
\begin{fait}[lemma 4,
 \cite{ana}]\label{amalgamationassy}
Let $A, B, C$ be L-structures, let $i$ be an immersion from $A$ into $B$, 
and  $h$ be a homomorphism from $A$ to $C$.
 There exist a model $D$ of $T_ k(C)$  
a homomorphism $h'$  from $B$ to $D$, 
and an immersion $i_m$ from $C$ into $D$ such that 
the following diagram commutes
\[
\xymatrix{
    A \ar[r]^{i} \ar[d]_{h} & {B} \ar[d]^{h'} \\
    C \ar[r]_{im} & {D}
  }
\]
\end{fait}
\begin{fait}[lemma 3, \cite{ana}]\label{caracmax}
Let $T$ be an $h$-inductive theory and $A$ a model of $T$.
Then the following properties are equivalent:
\begin{enumerate}
\item $A$ is an amalgamation basis,
\item For every $ a$ n-tuple from  A, there exists a unique type
in $S_n(T)$ that contains $F_A( a)$.
\end{enumerate}
\end{fait}

 \begin{fait}[lemma 9, \cite{ana}]\label{robamalg}
 Let $T$ be a positive Robinson theory. Then we have the following
 properties:
\begin{enumerate}
\item  Every model of T that embeds in a pc model of T is a h-maximal
model of $T$.
\item  The h-maximal models of T have the amalgamation property.
\end{enumerate}
Moreover, if $T$ is h-universal then these two conditions are sufficient to
conclude that $T$ is a positive Robinson theory.
 \end{fait}
 The following lemma is a generalization of the fact \ref{robamalg}. 
 It gives a connection between the  positive 
 Robinson theory and the amalgamation property.
 \begin{lem}
An h-inductive theory $T$ is a positive Robinson theory 
if and only if the class of substructures of pc models of $T$ has 
the amalgamation property.
\end{lem}
\preuve 
Suppose that $T$ is a positive Robinson theory, then 
 $T_u(T)$ is a positive Robinson Theory. According to lemma
\ref{robsubpec}, every substructures of pc models is a 
h-maximal model of $T_u(T)$;  by  the fact 
\ref{robamalg} the class of h-maximal models of 
$T_u(T)$ has the amalgamation property.

Suppose now that  the class of substructures of 
pc models of $T$ has the amalgamation property.
 Let $A$ be a substructure of a pc model $A_e$ of $T$.
We claim that $A$ is a h-maximal model of $T_u(T)$. Indeed, let  
$f$ be a homomorphism from $A$ into  a model  $B$ of 
$T_u(T)$, let  $B_e$ be a pc model of $T$ in which 
$B$ is continued  by a homomorphism $g$. 
Since the class of substructures of pc models 
has the amalgamation property, we obtain the following commutative
diagram 
 \[
\xymatrix{
    A \ar[rr]^{e} \ar[d]_{f}& & {A_e} \ar[d]^{i_1} \\
    B\ar[r]_{g}&  B_e\ar[r]_{i_2}& C
  }
\] 
Where $C$ is a substructure of a pc model, then  $C$ is  a 
model of $T_u(T)$. On the other hand, the fact that 
$A_e$ and $B_e$ 
are pc models entails that  $i_1, i_2$ are immersions. However,
as  
$i_2\circ g\circ f= i_1\circ e$,  and $i_1\circ e$ is 
an embeddings then 
the homomorphism $f$ is an embedding. Thereby $A$ is a 
h-maximal of $T_u(T)$.
 Therefore, the class of h-maximal models 
of $T_u(T)$ is the class of substructures of pc models of 
$T$. By the fact \ref{robamalg} $T_u(T) $ is a positive
 Robinson theory.
\begin{df}\label{def1}
Let $T$ be an h-inductive theory and  $p\in S(T)$. Let
$D$ be a set of 
positive formulas. $D$ is called $p$-dense
if and only if it satisfies the following property:
$$\forall q\in S(T),\, D\subseteq q \Rightarrow p=q.$$
\end{df}
\begin{rem}

 A theory $T$ is a positive Robinson theory if and only if for 
every type $p$ in $S(T)$ there is a set of quantifier-free positive formulas 
which is p-dense. 
\end{rem}

\begin{lem}\label{caracdens}
A model $A$ of A theory $T$ is an amalgamation basis if and only if for every
$ \bar a \in A$ there is 
$p\in S(T)$ such that  $ F_A(\bar  a)$ is p-dense. 
\end{lem}
\preuve 
Let $A$ be an amalgamation basis of $T$. By the fact \ref{caracmax}, 
for every $ \bar a\in A$ there is a unique type $p\in S_n(T)$ 
such that $p\models F_A(\bar  a)$, 
hence $ F_A(\bar a)$
is p-dense.

Suppose now that for any $\bar a \in A$ there is  
$p\in S_n(T)$ such that  $ F_A(\bar a)$ is p-dense, then $p$ is the 
unique type that contains $F_A( \bar a)$. By the fact \ref{caracmax},
$A$ is an amalgamation basis of $T$.

\begin{cor}
Let $A$ be an amalgamation bases and let $f$ be a homomorphism 
from $A$ into  a model $B$ of $T$.
Then  $f(A)$ is an amalgamation basis of $T_u(T)$.
\end{cor}

\begin{cor}
The class of amalgamation bases of an h-inductive theory 
is inductive.
\end{cor}
\preuve Let $I$ be a totally ordered set, 
and  let $(A_i, f_{ij})_{i,j\in I}$ be a inductive family of
  amalgamation bases of $T$. 
Let $A$ be the  inductive limit of $(A_i, f_{ij})_{i,j\in I}$.  
From the construction of $A$ we have, for every   $ \bar a\in A$
there is $i\in I$ such that $\bar  a\in A_i$ and 
$F_{A_i}(\bar a)\subset F_A(\bar a)$. Since $A_i$ 
is an amalgamation basis, there is
$p\in S(T)$ such that $F_{A_i}(\bar a)$ is p-dense.
Thereby $ F_{A}(\bar a)$ is p-dense, and by the lemma \ref{caracdens} 
$A$ is an amalgamation basis of $T$.

We denote by $T_{m}(T)$ the h-inductive theory of the class of 
amalgamation bases of $T$. We have 
$T\subset T_{m}(T)\subset T_k(T)$.\\
\begin{rem} 
\begin{itemize}
\item $T$ and $T_{m}(T)$ are companion theories.
\item The class of amalgamation bases is elementary if and only if it is axiomatized 
by $T_{m}(T)$ (fact \ref{ind}).
\end{itemize}
\end{rem}
\begin{thm}
A model $A$ of $T$ is an amalgamation basis if and only if
there is a pc model $A_e$ of $T$ in which $A$ is continued by a homomorphism 
$f$ and  such that the triple $(A, A_e, f)$ satisfies the following 
property;\\  for every 
positive formula $\varphi$ and $\bar a\in A$;  
if  $A_e\nvdash\varphi(f(\bar a))$ then 
there is $\varphi^c\in Res_T(\varphi)$ such that $A\models\varphi^c(\bar a)$.
\end{thm}
\preuve 
Suppose the property cited in the theorem
 holds for a model $A$ of $T$. We claim that
$A$ is an amalgamation basis of $T$. Indeed, let $\bar a\in A$. 
Pose $p=tp(f(\bar a))$ in $A_e$.
 We will  show that $F_a(\bar a)$ is 
p-dense.\\
 Let $q\in S(T)$ such that $F_A(\bar a)\subseteq q$. Let 
$\varphi$ be a positive formula such that $p\nvdash\varphi$, by hypothesis there is $\varphi^c$ a positive formula such that
$\varphi^c\in F_A(\bar a)$ and   
$T\vdash\neg \exists\bar x\,(\varphi(\bar x)\wedge\varphi^c(\bar x))$.
Then $q\vdash\varphi^c$ and  $q\nvdash\varphi$. Thereby 
$q\subseteq p$, by the maximality of types we obtain $p=q$.\\

For the reverse direction.  Let $A$ be an amalgamation 
 basis of $T$. Owing to 
 the lemma \ref{caracdens}, for every $\bar a\in A$ 
  there is a unique type $p\in S_n(T)$ such that $F_A(\bar a)$
 is p-dense. Let $\varphi$ be a positive formula such that 
 $p\nvdash\varphi$, then every  continuations of $A$ do not satisfy 
 $\varphi(\bar a)$ because  otherwise $A$ would not be an amalgamation
 basis.
 Thereby the set of h-inductive sentences  $\{T,  Diag^+(A),  \varphi(\bar a)\}$ is inconsistent. Then
 by compactness there exists $\psi(\bar a, \bar b)\in Diag^+(A)$
 such that $T\vdash \neg\exists\bar x\,(\varphi(\bar x)
 \wedge\exists\bar y\psi(\bar x, \bar y))$. Denote by $\varphi^c$ the 
 positive formula $\exists\bar y\psi(\bar x, \bar y)$. 
 we obtain $A\models\varphi^c(\bar a)$ and
$T\vdash \exists\bar x\,(\varphi(\bar x)\wedge\varphi^c(\bar x))$.
The following characterization of the Hausdorff  property of 
an h-inductive theory is given is \cite{begnacpoizat}.
 \begin{fait}[theorem 20, \cite{begnacpoizat}]\label{th20}
An h-inductive theory $T$ is Hausdorff if and only if 
the theory $T_k(T)$ has 
the amalgamation property.
\end{fait}
\begin{rem}
Considering $T_m(T)\subset _k(T)$, if the class of amalgamation bases is elementary then the theory
 is Hausdorff. 
\end{rem} 

\begin{cor}
If the class of h-maximal models of $T$ is elementary  and 
has the amalgamation property then $T$ is Hausdorff.
\end{cor}
\preuve 
Since $T\leq T^\star(T)\leq T_k(T)$ where $T^\star(T)$ is the h-inductive theory of 
the class of h-maximal models, thus
every models of $T_k(T)$ is a h-maximal model  and then an amalgamation basis of $T$. By the fact \ref{th20} $T$ is Hausdorff.

In the rest of this section we will  discuss  the preservation of the property
of amalgamation by extension and restriction.

\begin{lem}\label{resamalbasis}
Let $B$ be an amalgamation basis of an h-inductive theory $T$ and  
$A$ a structure  which is immersed in $B$, then $A$ is 
an amalgamation basis of $T$.
\end{lem}
\preuve
  Let $B_1, B_2$ be  models of 
 $T$,  $f$ and $g$ homomorphisms from $A$ respectively into
 $B_1$ and $B_2$.  Since $A$ is immersed in $B$, by the fact \ref{amalgamationassy},
 there exist $C_1\models T_k(B_1)$ and $C_2\models T_k(B_2) $ 
 ( which  are evidently  model of $T$) 
 such that the following diagram commutes.
 
\[
\xymatrix{
     & B_1\ar[r]^{f_1}& C_1\\
    A \ar[ru]^{f} \ar[rd]_{g}\ar[r]^{i_m} & {B}
     \ar[ru]^{h_1} \ar[rd]^{h_2} \\
    & B_2 \ar[r]_{g_1} & {C_2} 
  }
\]  
 where $h_1, h_2$ are  homomorphisms, $i_m, f_1$ and  $g_1$ 
 are immersions.\\ 
 Now, since $B$ in an amalgamation basis of $T$ and 
$C_1, C_2$ models of $T$, we obtain  the following commutative diagram
\[
\xymatrix{
     & B_1\ar[r]^{f_1}& C_1\ar[rd]^{f'}\\
    A \ar[ru]^{f} \ar[rd]_{g}\ar[r]^{i_m} & 
    {B} \ar[ru]^{h_1} \ar[rd]^{h_2}&  & D \\
    & B_2 \ar[r]_{g_1} & {C_2} \ar[ru]^{g'}
  }
\]    
Where  $f', g'$ are homomorphisms and $D\models T$. 
Thereby $A$ is an amalgamation basis of $T$.
\begin{fait}\label{monthm}(theorem 1, \cite{ana})
Let $A$ be a L-structure and  a $B$ pc model  of $T_k(A)$ 
(with parameters in $A$). $T_k(A)$ is Hausdorff 
if and only if $T_k(B)$ is Hausdorff.
\end{fait}
\begin{lem}
Let $T$ be a complete h-inductive theory and let $A_e$ and $B_e$ be  
pc models of $T$. $A_e$ is Hausdorff if and only if $B_e$ is Hausdorff. 
\end{lem}
\preuve 
Since $T$ is complete, there is  a pc model $C_e$ of $T$ such  
$A_e$ and $B_e$ are continued in $C_e$. 
On the other hand, as $A_e, B_e$  and $C_e$ are 
pc models of $T$, then  $C_e$ is respectively a 
pc model of $T_k(A_e)$  in the language obtained from the
 language of $T$ by adding
the elements of $A_e$ as constants;
  and  a pc model of  
$T_k(B_e)$ with parameters in $B_e$. By the fact 
\ref{monthm},  $B_e$ is Hausdorff. 
\begin{lem}
Let $T$ be a complete h-inductive theory. If  $T$ is Hausdorff
and  $A_e$ a pc model of $T$, then   $T_k(A_e)$
(with parameters in $A_e$) is Hausdorff.
\end{lem}
\preuve
Let $A_e$ be a pc model of $T$.
Suppose that $T$ is Hausdorff. We will show that $T_k(A_e)$ has 
the amalgamation property.\\
Let  
$B, C$ and $D$ be  models of  $T_k(A_e)$ such that $B$ is 
continued into $C$ and $D$ respectively by the homomorphisms 
$f$ and $g$. 
Since $B, C$ and $D$  are also 
 models of $T_k(T)$ and $T_k(T)$ has the amalgamation 
 property;  then there is  a model $E$ of 
 $T_k(T)$ such that the following diagram commutes;
 \[
\xymatrix{
    A_e\ar[r]^{i_m} & B\ar[r]^{f}\ar[d]_{g}& C\ar[d]^{f'}\\
     & D \ar[r]_{g'} & E 
  }
\] 
where $i_m$ is an immersion,  $f'$ and $g'$ are homomorphisms.\\
Now, since $E\vdash T$ and $A_e$ is a pc model of $T$, then 
$f'\circ f_{|A_e}$ and $g'\circ g_{|A_e}$ are immersions.
Using the fact that $A_e$ is immersed in $E$, 
 and $A_e$ is a pc model
of $T$ we obtain the consistency of the set 
of h-inductive sentences $\{T_k(A_e), Diag^+(E)\}$.
Indeed, in the contrary case by compactness we would have 
 a positive formula $\varphi$ and  $\bar a\in A_e$ such that 
$E\models\exists \bar x\ \varphi(\bar a, \bar x)$ and 
$T_k(A_e)\vdash\neg\exists\bar{x}\ \varphi(\bar a, \bar x)$, this 
entails $A_e\vdash\exists \bar x\ \varphi(\bar a, \bar x)$.
Contradiction with   $T_k(A_e)\vdash\neg\exists\bar{x}\ \varphi(\bar a, \bar x)$. Thus $\{T_k(A_e), Diag^+(E)\}$ is consistent, thereby 
$E$ is continued into $F$ a model of $T_k(A_e)$, and we obtain the 
commutative diagram 
$$
 \xymatrix{
      B\ar[r]^{f}\ar[d]_{g}& C\ar[d]^{f'}&\\
      D \ar[r]_{g'} & E\ar[r]^h& F 
  }
$$
Therefore $T_k(A_e)$ has the amalgamation property, and by the 
fact \ref{th20}, $A_e$ is Hausdorff. 
\section{Characterization  of the amalgamation basis and amalgamation property}
\begin{thm}\label{caracbaseamalgamationaspec}
Let $A$ be a model of an h-inductive theory $T$. 
$A$ is 
an amalgamation basis of $T$ if and only if it satisfies the 
following property:\\
For every positive formula $\phi$ and $\bar{a}$ a tuple 
of $A$; if $A\nvDash\phi(\bar a)$ then there exist
$\phi'\in Res_T(\phi)$, and 
$(\varphi_\psi)_{\psi\in Res_T(\phi)}$ a family 
of positive formulas such that; for every 
$\psi\in Res_T(\phi)$ we have $\varphi_\psi\in Res_T(\psi)$, and: 
$$A\models \phi'(\bar{a})\vee \bigwedge_{\psi\in Res_T(\phi)}
\varphi_\psi(\bar{a}).$$
\end{thm}
\preuve 
Let $A$ be an amalgamation basis of $T$ and let $\phi$ be a positive formula. Consider $\bar a$  a tuple of $A$ such that $A\nvDash\phi(\bar a)$. Then
we distinguish two cases depending on whether the set of h-inductive sentences $\{T, Diag^+(A), \phi(\bar a)\}$ is consistent or not.
 \begin{itemize}
\item $\{T, Diag^+(A), \phi(\bar a)\}$ is inconsistent.  By compactness 
there is  a positive formula $\phi'$ such that 
$A\vDash\phi'(\bar a)$ and 
$T\vdash\neg\exists\bar x\ (\phi(\bar x)\wedge\phi'(\bar x))$.
Thereby $\phi'\in Res_T(\phi)$ and $A\vDash\phi'(\bar a)$.
\item $\{T, Diag^+(A), \phi(\bar a)\}$ is 
consistent.
Since $A$ is an amalgamation basis, then
  for every $\psi\in Res_T(\phi)$ the set 
of h-inductive sentences $\{T, Diag^+(A), \psi(\bar a)\}$ 
is inconsistent.
By compactness there is $\varphi_\psi\in Res_T(\psi)$
such that $A\models\varphi_\psi(\bar a)$. Hence we obtain
$A\vDash\bigwedge_{\psi\in Res_T(\phi)}\varphi_\psi(\bar a)$.
\end{itemize}

For the reverse direction, consider $A$  a model of $T$ such that; 
For every positive formula $\phi$ and $\bar{a}$ a tuple 
of $A$, if $A\nvDash\phi(\bar a)$ then there exist 
$\phi'\in Res_T(\phi)$, and  
$(\varphi_\psi)_{\psi\in Res_T(\phi)}$ a family 
of positive formulas such that. For every 
$\psi\in Res_T(\phi)$ we have $\varphi_\psi\in Res_T(\psi)$, and: 
$$A\vdash \phi'(\bar{a})\vee \bigwedge_{\psi\in Res_T(\phi)}
\varphi_\psi(\bar{a})$$
Suppose that $A$ is not an amalgamation basis of $T$, then 
there exist $B$ and $C$ pc models  of $T$ in which $A$ is 
continued and such that,  in 
the language obtained from the language of $T$ by adding
the elements of $A$ as constants; the set of 
h-inductive sentences
$\{T, Diag^+(B), Diag^+(C)\}$
is inconsistent. By compactness
there are $(\phi, \psi)$ a pair of positive formulas and $\bar a\in A$ such that $\varphi\in Res_T(\phi)$ and;
$$B\models\phi(\bar a),\ \ \  C\models\varphi(\bar a),\ \ \ 
 A\nvDash\phi(\bar a)$$
By hypothesis, there exist $\phi'\in Res_T(\phi)$ and 
 $(\varphi_\psi)_{\psi\in Res_T(\phi)}$ a family 
of positive formulas such that; for every 
$\psi\in Res_T(\phi)$ we have $\varphi_\psi\in Res_T(\psi)$, and: 
$$A\vdash \phi'(\bar{a})\vee \bigwedge_{\psi\in Res_T(\phi)}
\varphi_\psi(\bar{a})$$
Suppose that  $A\models\phi'(\bar{a})$, then 
$B\models \phi'(\bar{a})\wedge\phi(\bar{a})$ which is impossible. Consequently 
$A\models \bigwedge_{\psi\in Res_T(\phi)}
\varphi_\psi(\bar{a})$. Moreover, 
since $A$ is continued into $C$,  then 
$C\models  \bigwedge_{\psi\in Res_T(\phi)}
\varphi_\psi(\bar{a})$. By the fact 
that $C$ is a pc model of $T$, and for every $\psi\in Res_T(\phi)$
we have $C\models\varphi_\psi(\bar{a})$  then $C\models \phi(\bar a)$.
Contradiction with  $C\vdash\varphi(\bar a)$.
Therefore $A$ is an amalgamation basis of $T$.
\begin{thm}\label{th21}
The class of amalgamation bases of an h-inductive theory $T$ 
is elementary if and only if it is axiomatized 
by  the following  h-inductive theory:\\ 
For every positive formula $\phi$
and for every $\phi'\in Res_T(\phi)$ there is
$(\varphi, \psi)$ 
a pair of $T_m(T)$-weakly complementary positive formulas such that
 $\varphi\in Res_T(\phi)$ and $\psi\in Res_T(\phi')$.
\end{thm}
\preuve
Suppose that the class of amalgamation bases is elementary and assume the existence of a pair of positive formulas $(\phi, \phi')$ 
such that $\phi\in Res_T(\phi')$, and for every pair of 
positive formulas $(\varphi, \psi)$; if $\varphi\in Res_T(\phi)$ 
and $\psi\in Res_T(\phi')$ then 
$T_m(T)\nvdash \forall\bar x(\varphi(\bar x))\vee\psi(x))$.\\
By compactness we obtain the consistence of the set 
of the h-inductive sentences 
$$T^\star=
\{T_m(T), \neg\phi(\bar x), \neg\phi'(\bar x), \neg\varphi(\bar x), 
\neg\psi(\bar x)\ |\ \varphi\in Res_T(\phi); \psi\in Res_T(\phi')\}$$
Now, let $A$ be a model of $T^\star$, and $\bar a$ a tuple of 
$A$ such that
 $$A\nvDash\phi(\bar a)\vee\phi'(\bar a)\vee\varphi(\bar a)
 \vee\psi(\bar a)$$
 Where $\varphi$ and $\psi$ range  respectively through  
 $ Res_T(\phi)$ and $Res_T(\phi')$.\\
 By compactness we get the following diagram 
$$
\xymatrix{
A\ar[r]^{f}\ar[d]_{g}& B\\
C&
}
$$ 
where $f$ and $g$ are homomorphisms,  
$B$ and $C$ are models of $T$ such that $B\models\phi(f(\bar a))$ and   $C\models \phi'(\bar a)$.  It is obvious that the diagram thus obtained can not be completed with a model of $T$ to have
 a commutative diagram. This leads us to a contradiction.
  Thereby for every  pair of 
positive formulas $(\phi, \phi')$ such that 
$\phi'\in Res_T(\phi)$ there is
$(\varphi, \psi)$ 
a pair of $T_m(T)$-weakly complementary positive formulas such that
 $\varphi\in Res_T(\phi)$ and $\psi\in Res_T(\phi')$.

 For the reverse direction, let $A$ be a model of $T_m(T)$ 
 (the h-inductive theory of the class of amalgamation bases of $T$)
such that for every a pair of 
positive formulas $(\phi, \phi')$ with 
$\phi'\in Res_T(\phi)$, there is
$(\varphi, \psi)$ 
a pair of $T_m(T)$-weakly complementary positive formulas such that
$\varphi\in Res_T(\phi)$ and $\psi\in Res_T(\phi')$.\\
Suppose that $A$ is not an amalgamation basis of $T$. Then there exist
$B$ and $C$  model  of $T$ in which $A$ is continued, 
 $\bar a\in A$, $\phi(\bar x)$, and  $\phi'(\bar x)$ positive 
 formulas with 
$\phi'\in Res_T(\phi)$, such that $B_e\models\phi(\bar a )$
and $C_e\models\phi'(\bar a)$. Let $(\varphi, \psi)$ be a 
pair of $T_m(T)$-weakly complementary formula such that 
$\varphi\in Res_T(\phi)$ and $\psi\in Res_T(\phi')$. Suppose 
that $A\models\varphi(\bar a)$, thereby 
$B\models\phi(\bar a)\wedge\varphi(\bar a)$. Contradiction
with  $\varphi\in Res_T(\phi)$. The same is if 
$A\models\psi(\bar a)$. Therefore $A$ is an amalgamation basis
of $T$.
\begin{cor}
Let $T$ be a h-inductive theory. 
$T_k(T)$ has the amalgamation property if and only if it satisfies 
the following property.\\
For every positive 
formula $\phi$, and for every  $\phi'\in Res_T(\phi)$;
there exists a pair $(\varphi, \psi)$ of $T_k(T)$-weakly 
complementary positive formulas such that; 
$\varphi\in Res_T(\phi')$, and $\psi\in Res_T(\phi)$.
\end{cor}
\preuve
From the fact \ref{th20} and the theorem \ref{th21}.

\bibliographystyle{plain}

Mohammed  Belkasmi\\
Math  Department\\
College of Science\\
Qassim University

\end{document}